\input amstex
\documentstyle{amsppt}
\nopagenumbers \centerline{\bf Inverse Scattering for Gratings and
Wave Guides} \vskip.2in \centerline{Gregory Eskin, Department of
Mathematics, UCLA, E-mail: eskin\@math.ucla.edu} \vskip.1in
\centerline{James Ralston, Department of Mathematics, UCLA,
E-mail: ralston\@math.ucla.edu} \vskip.1in
 \centerline{Masahiro Yamamoto,
Graduate School of Mathematical Sciences,}

\centerline{University of Tokyo, E-mail: myama\@ms.u-tokyo.ac.jp}

\vskip.2in \noindent{\bf Abstract:} We consider the problem of
unique identification of dielectric coefficients for gratings and
sound speeds for wave guides from scattering data. We prove that
the \lq\lq propagating modes" given for all frequencies uniquely
determine these coefficients. The gratings may contain conductors
as well as dielectrics and the boundaries of the conductors are
also determined by the propagating modes. \vskip.2in
 \noindent \S 0. {\bf Introduction}
\vskip.1in
 Consider Maxwell's equations
for time-harmonic electric and magnetic fields, $\exp(-i\omega
t)E(x_1,x_2,x_3)$ and $\exp(-i\omega t)H(x_1,x_2,x_3)$, in the
absence of currents and charges
$$\nabla \times E-i\omega\mu_0 H=0,$$
$$\nabla \times H+i\omega \epsilon(x) E=0.$$
In this paper we study the inverse problem of determining the
electric permittivity, $\epsilon$, and hence the dielectric
coefficient, $\epsilon/\epsilon_0$, from scattering data for these
equations. The fundamental assumptions are that $\epsilon$ is
independent of $x_3$, $2\pi$-periodic in $x_1$ and  constant
($=\epsilon_0$) for $|x_2|> T.$ These conditions are designed to
model a dielectric \lq\lq grating" extending throughout the region
$|x_2|< T$ (c.f. [P], [BDC]). We also allow for conducting bodies
embedded in the dielectric as long as they satisfy conditions
analogous to our conditions on $\epsilon$: their boundaries should
be invariant with respect to all translations in $x_3$ and
translations by $2\pi$ in $x_1$, and they should be contained in
$|x_2|< T$.

To define data sets for this inverse problem it is customary to
consider the scattering problem for fields with either the
transverse electric (TE) or transverse magnetic (TM)
polarizations, respectively $E(x)=(0,0,u(x_1,x_2))$ and
$H(x)=(0,0,v(x_1,x_2))$. These polarizations reduce Maxwell's
equations, respectively, to
$${1\over \epsilon(x)}\Delta u +\omega^2\mu_0u=0 \eqno{(TE) }$$
and $$\nabla \cdot {1\over \epsilon(x)}\nabla
v+\omega^2\mu_0v=0.\eqno{(TM)}$$ In the case of embedded
conductors we consider the TE polarization in the exterior of the
conductors with the Dirichlet condition, $u=0$, on the boundary,
since this corresponds to $E=0$ and $H\cdot\hat n=0$ on the
surfaces of the conductors. For our purposes it is convenient to
write both (TE) and (TM) as
$$Lu-k^2u=0, $$
where $L=-\Delta$ for $|x_2|>T$ and $k^2=\omega^2\mu_0\epsilon_0$.

  We will also present the analogous inverse problem for
acoustic wave guides. This requires only small modifications of
the arguments for gratings. The wave guides that we consider are
simply slabs, $\{0<x_1<B\}$, in which the sound speed $c$ is a
function of $(x_1,x_2)$. We assume that $c(x_1,x_2)=c_0(x_1)$ for
$|x_2|>T$, and impose Dirichlet condition on $x_1=0$, and the
Neumann condition, $\partial_{x_1}u=0$, on $x_1=B$. These boundary
conditions correspond to an acoustically soft reflecting surface
at $x_1=0$ and an acoustically hard reflecting surface at $x_1=B$,
modelling underwater sound propagation with $x_1$ as depth (c.f.
[BGWX]). We will show that scattering data from propagating modes
for the operator $L=-c^{2}(x)\Delta$ with these boundary
conditions determine $c(x)$.

In both these settings we will apply recent results on inverse
coefficient problems for hyperbolic equations (Belishev [B],
Kachalov-Kurylev-Lassas [KKL] and Eskin [E1],[E2]). In those
papers the data for the inverse problem is the
Dirichlet-to-Neumann map. Hence the objective here will be to show
that the scattering data determine the Dirichlet-to-Neumann map on
a line $x_2=T$.

Inverse scattering problems for dielectric gratings have been
studied previously in [BDC], [BF], [K], [HK] and [EY]. These
articles consider primarily the inverse problem of finding the
boundaries of conductors embedded in a dielectric of constant
permitivity from scattering data. To the best of our knowledge the
present paper is the first to show that a variable dielectric is
uniquely determined by scattering data.

Inverse coefficient problems for wave guides were studied in
[BGMX], [GMX], [M], [X] and [DM].  These papers give methods for
recovering the sound speed. We only consider the uniqueness
problem and prove that the sound speed is uniquely determined by
the propagating modes. Our approach was influenced by the work of
S.Dediu and J. McLaughlin, [DM], which also uses propagating
modes.
 \vskip.2in \noindent \S 1. {\bf
Statement of results}
 \vskip.1in
Our results for gratings hold under mild conditions on the
operator on $L$ in (1). We will assume that $L$ is a second order
elliptic operator on $D\subset {\Bbb R}^2$ which is symmetric in
the inner product
$$ (f,g)=\int_Df(x)\overline {g(x)} a(x)dx.$$ The weight $a(x)$ is
continuous and strictly positive on $\overline D$. The
coefficients of $L+\Delta$ are supported in $|x_2|< T$, and $L$
commutes with translation by $2\pi$ in $x_1$. Likewise $a(x)-1$ is
supported in $|x_2|<T$ and $a(x_1+2\pi,x_2)=a(x_1,x_2)$. We will
also assume that the region $D$ is invariant under translation by
$2\pi$ in $x_1$, and boundary $D$ is smooth. There are two cases
that we wish to consider. \vskip.1in \noindent Case 1: $D$ is
connected and contains $\{|x_2|>T\}$. In other words, while there
may be some holes in $D$, they do not disconnect $D$, and they are
contained in $|x_2|<T$. \vskip.1in \noindent Case 2: $D$ is
connected and we have the inclusions
$$\{x_2>T\}\subset D\subset \{x_2>-R\}$$
for some $R>0$.

The domain of $L$ will be $H^2_\alpha(D)\cap H^1_{\alpha,0}(D)$.
By $H^k_{\alpha}(D)$ (resp. $H^k_{\alpha,0}(D)$) we mean functions
satisfying
$$u(x_1+2\pi,x_2)=e^{2\pi i\alpha}u(x_1,x_2)\eqno{(1)}$$
such that $\phi(x_1)u(x_1,x_2)\in H^k(D)$ (resp. $H^k_0(D)$) for
all $\phi\in C_c^\infty(\Bbb R)$. Note that this domain for $L$
corresponds to the Dirichlet boundary condition on $\partial D$.

For wave guides we simply take $L=-c(x)^{2}\Delta$ on
 $D=\{x\in {\Bbb R}^2:0<x_1<B\}$, with
Dirichlet and Neumann boundary conditions on $x_1=0$ and $x_1=B$,
respectively. As indicated above, $c(x)=c_0(x_1)$ when $|x_2|>T$.
\vskip.1in

For both gratings and wave guides scattering data at fixed energy
$k^2$ are obtained from the \lq\lq propagating modes". In the case
of gratings the scattered wave $v_+$ is obtained by solving
$(L-k^2)u=0$ in $D$ with $u=\exp (l\cdot x)+v_+$, where $l\cdot
l=k^2$ and $v_+$ is the \lq\lq outgoing" solution to
$$(L -k^2)v=-(L-k^2)e^{il\cdot x}$$
obtained as the limit as Im$\{k\}\to 0_+$ (see below). For
gratings we will only use scattering data from incident waves
$\exp(il\cdot x)$ which satisfy the condition (2), i.e. those with
$l_1=n+\alpha, n\in \Bbb Z,\ l_2=-\sqrt{k^2-(m+\alpha)^2} $. For
$x_2>T$ the scattered wave $v_+$ has the form
$$v_+(x,n,k)=\sum_{m\in {\Bbb
Z}}e^{i[(m+\alpha)x_1+x_2\sqrt{k^2-(m+\alpha)^2}]}a_m(n,k)\eqno{(2)}$$
$$=
\sum_{(m+\alpha)^2<k^2}e^{i[(m+\alpha)x_1+x_2\sqrt{k^2-(m+\alpha)^2}]}a_m(n,k)+\hbox{O}(e^{-\delta
x_2})$$ for some $\delta>0$, provided that $k$ does not belong to
the set of \lq\lq thresholds ", $\{ k: k^2=(p+\alpha)^2, p\in
{\Bbb Z}\}$. We call $\{a_m(n,k),\ n,m\in {\Bbb Z}:\
(m+\alpha)^2<k^2, (n+\alpha)^2<k^2\}$ the scattering data at
energy $k^2$ from \lq\lq propagating modes".

For wave guides, since we are taking $L_0=-c_0(x_1)^{2}\Delta$ as
the unperturbed operator, the scattered wave $v_+(x,k,n)$  is
obtained by solving $(L-k^2)u=0$ with $u=\Phi(x,k,n)+v_+$, where
$\Phi$ is a generalized eigenfunction for $L_0$, i.e.
$\Phi(x,k,n)=\exp(-ix_2\sqrt{\mu_n(k)})\phi_n(x_1,k)$, where
$$\phi_n^{\prime\prime} +k^2c_0^{-2}\phi_n=\mu_n(k)\phi_n,\
\phi_n(0)=0,\ \phi_n^\prime(B)=0,\hbox{ and }\mu_n(k)>0.$$  For
$x_2>T$ the scattered wave $v_+$ has the form
$$v_+(x,k,n)=\sum_{m\in \Bbb
N}b_m(k,n)e^{ix_2\sqrt{\mu_m(k)}}\phi_m(x_1,k)$$
$$\sum_{\mu_m(k)>0}b_m(k,n)e^{ix_2\sqrt{\mu_m(k)}}\phi_m(x_1,k)+\hbox{O}(e^{-\delta
x_2})\eqno{(2^\prime)}$$ for some $\delta>0$, provided that $k$
does not belong to the set of thresholds, $\{k: \mu_p(k)=0, p\in
{\Bbb N}\}$. Here we call $\{b_m(k,n),\ n, m\in{\Bbb N}:\
\mu_n(k)>0,\ \mu_m(k)>0\}$ the scattering data at energy $k^2$
from propagating modes.

Note that in these definitions the functions
$\sqrt{k^2-(m+\alpha)^2}$ and $\sqrt{\mu_m(k)}$ will be chosen so
that they extend into Im$\{k\}>0$ with positive imaginary parts.
Letting $z(k)$ stand for either $k^2-(m+\alpha)^2$ or $\mu_m(k)$,
this choice amounts to choosing $\sqrt{z(k)}>0$ when $z(k)>0$ and
$k>0$, $\sqrt{z(k)}<0$ when $z(k)>0$ and $k<0$, and
$\sqrt{z(k)}=i\sqrt{|z(k)|}$ when $z(k)<0$. We will follow these
conventions in the rest of the paper.

With the preceding definitions we have:
 \vskip.1in \noindent {\bf Theorem 1:} The scattering data from
 propagating modes in $x_2>T$ given for all $k$ determine $D$ and $\epsilon(x)$
 for gratings with either the (TE) or (TM) polarizations, and $c(x)$
 for wave guides.
 \vskip.1in
The proof of theorem will proceed as follows. We will consider
\lq\lq generalized distorted plane waves"
$$u_+(x,k,n)=\exp(i(n+\alpha)x_1-ix_2\sqrt{k^2-(n+\alpha)^2}) 
+v_+(x,k,n)$$ and
$$u_+(x,k,n)=\exp(-ix_2\sqrt{\mu_n(k)})\phi_n(x,k) +v_+(x,k,n),$$ which 
are
defined without the restrictions $k^2>(n+\alpha)^2$ and
$\mu_n(k)>0$. These generalized distorted plane waves exist for
$k\in \Bbb R\backslash S$, where $S$ is a discrete set. Note that
when $k^2<(n+\alpha)^2$ or $\mu_n(k)<0$ these generalized
distorted plane waves grow exponentially as $x_2\to\infty$. In \S
2 and \S 3 we show that the set of generalized distorted plane
waves, given for a fixed $k$ and all $n$, uniquely determines the
Dirichlet-to-Neumann map on a suitable line $x_2=T$ for all
choices of  $k$ outside a discrete set. We also show that, making
use of the analytic continuation to Im$\{k\}>0$ of the
$v_+(x,k,n)$'s,  these generalized distorted plane waves are
determined by the scattering data from propagating modes. Thus,
under the hypotheses of Theorem 1, the Dirichlet-to-Neumann map is
known on $x_2=T$ for all $k\in \Bbb R$ outside a discrete set.
Since this is equivalent to knowing the {\it hyperbolic}
Dirichlet-to-Neumann map for the wave equations $u_{tt}=Lu$ on
$x_2=T$, the proof of Theorem 1 will be reduced to the results on
hyperbolic inverse coefficient problems cited above. Since
analytic continuation plays a big role here, there are many
variations on the set of $k$ for which the propagating modes are
known which lead to the same results. \vskip.2in

\noindent \S 2. {\bf Determination of the Dirichlet-to-Neumann Map
for Gratings} \vskip.1in
 In this section we will show that the scattering data from
 propagating modes
 determine the Dirichlet-to-Neumann map on a line $x_2=T$ for the case 
of gratings.
To do this we will first show that the traces of an appropriate
family of distorted plane waves on $x_2=T$ are dense in
$L^2(0<x_1<2\pi)$.

To begin we need the incoming and outgoing fundamental solutions
for $-\Delta-k^2$ on ${\Bbb R}^2$ in a form compatible with (1).
Using Fourier series in $x_1$ to reduce this to an ODE in $x_2$,
one computes that for Im$\{k\}>0$
$$[(-\Delta -k^2)^{-1}f](x)=$$
$$\sum_{m=-\infty}^\infty {1\over 2\pi}\int_{-\infty}^\infty
(\int_0^{2\pi}e^{i[(m+\alpha)(x_1-y_1)+\lambda_m(k) |x_2-y_2|]}{
f(y)\over 2i\lambda_m(k)}dy_1)dy_2,\eqno{(3)}$$ where
$$\lambda_m(k)=\sqrt{k^2-(m+\alpha)^2}$$ with the
branch chosen  $\sim k$ near infinity and the cut on
$(-|m+\alpha|,|m+\alpha|)$. Note that \noindent
Im$\{\lambda_m(k)\}>0$ for Im$\{k\}>0$, and hence $(-\Delta
-k^2)^{-1}$ maps $H^0_\alpha({\Bbb R}^2)$ into $H^2_\alpha({\Bbb
R}^2)$. The continuous extension of $\lambda_m(k)$ from
Im$\{k\}>0$ to the real axis is given by
$$\lambda_m(k)=k\sqrt{1-(m+\alpha)^2/k^2}$$ when $(m+\alpha)^2<k^2$
and $$\lambda_m(k)=i\sqrt{(m+\alpha)^2-k^2}$$ when
$(m+\alpha)^2>k^2$. The corresponding extension of
$(-\Delta-k^2)^{-1}$ to ${\Bbb R}\backslash\{\pm(m+\alpha),\ m\in
{\Bbb Z}\}$ gives the outgoing fundamental solution, $G_+(k)$. For
the incoming fundamental solution we take
$$\lambda_m(k)=-\sqrt{k^2-(m+\alpha)^2},$$
i.e. the branch chosen  $\sim -k$ near infinity. Substituting this
in the formula for $(-\Delta -k^2)^{-1}f$ to get
$(-\Delta-k^2)^{-1}$ in Im$\{k\}<0$, and define the incoming
fundamental solution, $G_-(k)$, by continuous extension from
Im$\{k\}<0$ to the real axis. Hence, by construction $G_+(k)$
extends analytically to $(-\Delta -k^2)^{-1}$ in Im$\{k\}>0$, and
$G_-(k)$ extends analytically to $(-\Delta -k^2)^{-1}$ in
Im$\{k\}<0$.

Now we turn to the construction of generalized distorted plane
waves for $L$. Choose $\psi\in C^\infty(\Bbb R)$ such that $\psi
\equiv 1$ on a neighborhood of $\{|x_2|\geq T\}$, and the support
of $\psi(x_2)$ is contained in the set where $L=-\Delta$ and
$a=1$. An (outgoing) generalized distorted plane wave for $L$ is a
solution of $(L-k^2)u=0$ in $D$ of the form
$u_+=\psi(x_2)\exp(il\cdot x)+v_+$, where $l_1\equiv \alpha$ mod 1
with $0\leq \alpha< 1,\ l_2=-\sqrt{k^2-(n+\alpha)^2}$, and $v_+$
is defined by limiting absorption, i.e.
$$v_+=-\lim_{\epsilon\to
0_+}(L-(k+i\epsilon)^2)^{-1}(L-k^2)(\psi(x_2)e^{il\cdot x}).$$
These are {\it generalized} distorted plane waves in the sense
used in \S 1, since the second component of $l$ is not necessarily
real.  Ordinarily,  outgoing distorted plane waves are defined as
solutions of the form $u_+=\exp(il\cdot x) + v_+$ where $v_+$ is
outgoing. However, we have
$$u_+=\psi(x_2)e^{il\cdot x}+v_+$$
$$=e^{il\cdot x}+[(\psi(x_2)-1)e^{il\cdot x} -\lim_{\epsilon\to
0_+}(L-(k+i\epsilon)^2)^{-1}(L-k^2)(\psi(x_2)e^{il\cdot x})].$$
Since $(\psi(x_2)-1)\exp(il\cdot x)$ is  outgoing, the term in
brackets is outgoing.

 The limit defining $v_+$ will
exist unless \vskip.1in \noindent i) $k$ is one of the \lq\lq
thresholds", $k^2=(n+\alpha)^2$, where $G_\pm(k)$ are undefined,
or \vskip.1in \noindent ii) there is a solution to the homogeneous
equation $(L-k^2)u=0$ in $D$ which is square-integrable on
$D\cap\{0<x_1<2\pi\}$.  \vskip.1in \noindent We denote the set of
exceptional $k$'s defined by i) and ii) as $S$.\footnote{Note that
case ii) can occur. Choose $V$ with compact support so that
$-\partial^2_{x_2}+V(x_2)$, considered as a Schr\"odinger operator
on $\Bbb R$, has a bound state, $u\in L^2(\Bbb R)$, i.e.
$(-\partial^2_{x_2}+V(x_2))u=Eu.$ Then, taking $m$ large enough
that $(m+\alpha)^2+E-V$ is strictly positive, defining
$\epsilon^{-1}(x)= ((m+\alpha)^2+E)((m+\alpha)^2+E-V)^{-1}$ and
$\psi=\exp(i(m+\alpha)x_1)u(x_2)$, we have $(1/\epsilon(x))\Delta
\psi+((m+\alpha)^2+E)\psi=0.$}

 Since $l_1=n+\alpha$ for a unique
$n\in \Bbb Z$, we use $n$ and $k$ to parametrize the generalized
distorted plane waves, $u=u(x,n,k)$. With these definitions we
have outgoing distorted plane waves for all $(n,k)\in \Bbb Z\times
\Bbb R\backslash S$. The analytic properties of $G_+(k)$ discussed
above carry over to the $u_+(x,n,k)$'s: they have analytic
continuations to Im$\{k\}>0$ which extend continuously back to
${\Bbb R}\backslash S$. This leads directly to the following
conclusion which we state as a lemma. \vskip.1in \noindent {\bf
Lemma 0:} For each $n$ the set $\{u_+(x,n,k),\ k\in I\}$, where
$I$ is an open interval in $k^2>(n+\alpha)^2$ determines
$u_+(x,n,k)$ for $k\in {\Bbb R}\backslash S$. Thus the true
distorted plane waves determine the generalized distorted plane
waves. \vskip.2in

The following observation is the main step in the proof.
\vskip.1in \noindent {\bf Lemma 1.}  Letting $(L-(k-i0)^2)^{-1}g$
denote $\lim_{\epsilon\to 0_-} (L-(k+i\epsilon)^2)^{-1}g$, the
\lq\lq incoming" solution, we have
$$\int_0^{2\pi}\overline {f(x_1)}u_+(x_1,T,m,k)dx_1=$$
$$\int_{D\cap \{0<x_1<2\pi\}}e^{i[(m+\alpha)x_1-x_2\lambda_m(k)]}
\overline{(L-k^2)\psi(L-(k-i0)^2)^{-1}(f\delta_T)}dx,\eqno{(4)}$$
where $\delta_T(\phi)=\int_0^{2\pi} \phi(x_1,T)dx_1,\
f(x_1)\in L^2(0<x_1<2\pi)$. \vskip.1in \noindent
Proof: We have
$$u_+=\psi e^{il\cdot x}-(L-(k+i0)^2)^{-1}(L-k^2)(\psi e^{il\cdot
x}).$$ Hence, letting $D_0=D\cap\{(x_1,x_2): 0<x_1<2\pi\}$, for
any smooth $g$ satisfying (1) with bounded support in $x_2$
$$\int_{D_0}u_+\overline g dx=\int_{D_0}\psi
e^{il\cdot x}\overline g dx -\int_{D_0}\overline 
g(L-(k+i0)^2)^{-1}(L-k^2)(\psi e^{il\cdot
x})dx.\eqno{(5)}$$
Since
$(L-(k-i0)^2)^{-1}$
is the adjoint of $(L-(k+i0)^2)^{-1}$, we have
$$\int_{D_0}\overline g(L-(k+i0)^2)^{-1}(L-k^2)(\psi e^{il\cdot
x})dx=\int_{D_0}\overline {(L-(k-i0)^2)^{-1}g}(L-k^2)(\psi e^{il\cdot
x})dx.\eqno{(6)}$$
Since $L=-\Delta$ on the support of $\psi$, for any  smooth $h$  
satisfying (1)
$$\int_{D_0}\overline h(L-k^2)(\psi e^{il\cdot
x})dx=\int_{D_0}e^{il\cdot x}\overline{
(2\nabla\psi\cdot\nabla+\Delta \psi)h}dx.\eqno{(7)}$$
Beginning with (5) and using (6) and (7), we have
$$\int_{D_0}u_+\overline g dx=\int_{D_0}\psi
e^{il\cdot x}\overline g dx-\int_{D_0}\overline 
{(L-(k-i0)^2)^{-1}g}(L-k^2)(\psi e^{il\cdot x})dx=$$
$$\int_{D_0}\psi
e^{il\cdot x}\overline g dx+\int_{D_0}\overline {(\nabla 
\psi\cdot\nabla +\Delta \psi)(L-(k-i0)^2)^{-1}g}e^{il\cdot x}dx=$$
$$\int_{D_0}
e^{il\cdot x}\overline{(L-k^2)\psi (L-(k-i0)^2)^{-1}g}dx$$ Now
approximating $f(x_1)\delta_T$ by $g$ of the form above gives (4).
\vskip.1in With (4) we can easily prove \vskip.1in \noindent {\bf
Lemma 2}: Assume that $k\in{\Bbb R}\backslash (S\cup S_T)$ is
fixed, where $S_T$ is the set of $k$ for which there are
nontrivial solutions to $Lu-k^2u=0$ which vanish on $x_2=T$ and
are square-integrable on $D\cap\{x_2<T\}$. Then the linear span of
$\{u_+(x_1,T,m,k),m\in \Bbb Z\}$ is dense in $L^2(0<x_1<2\pi)$.
 \vskip.1in \noindent Proof: Suppose that $f\in
L^2(0<x_1<2\pi)$ is orthogonal to the span of

\noindent $\{u_+(x_1,T,m,k),m\in \Bbb Z\}$. Then (4) implies
$$\int_{-\infty}^\infty(\int_0^{2\pi}e^{i[(m+\alpha)x_1-x_2\lambda_m(k)]}
\overline{(-\Delta-k^2)\psi(L-(k-i0)^2)^{-1}(f\delta_T)}dx_1)dx_2=0,$$
for all $m \in \Bbb Z$. Let
$$w=_{\hbox{def}}(L-(k-i0)^2)^{-1}(f\delta_T).$$
Since $w$ is incoming, we have $\psi w=G_-(k)(L-k^2)\psi w$.
Moreover, we have
$$(-\Delta-k^2)\psi w=f\delta_T-
(2\nabla\psi\cdot\nabla+\Delta \psi)(L-(k-i0)^2)^{-1}f\delta_T=0$$
for $x_2>T$. So when we represent $\psi w$ as $G_-(k)(L-k^2)\psi
w$ using the analog of (3) for $G_-(k)$, the integrand is
supported in $y_2\leq T$. Therefore, when $x_2> T$,
$|x_2-y_2|=x_2-y_2$ on the support of the integrand in (3), and
the identity above implies
$$\psi w(x)=0$$
for $x_2>T$, i.e. $w(x)=0$ for $x_2>T$.

At this point the arguments for Case 1 and Case 2 separate. In
Case 1, $w$ is an incoming solution to the homogeneous problem
$(L-k^2)w=0$ on $D\cap\{x_2<T\}$, satisfying (1) and $w(x_1,T)=0$.
Thus $w=0$ on all of $\partial (D\cap\{x_2<T\})$. In this case we
have for $R$ sufficiently large
$$0=\int_{D\cap\{0<x_1<2\pi\}\cap \{-R<|x_2|<T\}}
((\overline{Lw-k^2w})w-\overline w(Lw-k^2w))a(x)dx=$$
$$
\int_0^{2\pi}w{\partial \overline w\over\partial x_2}- \overline
w{\partial w\over\partial x_2}dx_1|_{x_2=-R}=\sum_{
(m+\alpha)^2<k^2}2\pi i\sqrt{k^2-(m+\alpha)^2}|a_m|^2
+O(e^{-\delta R}).\eqno{(8)}$$ The last equality comes from the
representation of  $\psi w$ as  $G_-(k)((L-k^2)\psi w)$, i.e. for
$x_2<-T$
$$w=\psi w=\sum_{\{m:
(m+\alpha)^2<k^2\}}e^{i(m+\alpha)x_1
-ix_2\sqrt{k^2-(m+\alpha)^2}}a_m +O(e^{-\delta |x_2|}).$$
 From
(8) it follows that the coefficients $a_m$ of the propagating
modes in $w$ vanish, and   $w\in L^2(D\cap\{0<x_1<2\pi\}\cap
\{x_2<T\})$. In other words $w$ is a Dirichlet eigenfunction for
$L$ in $D\cap \{x_2<T\}$ with the periodicity condition (1).

In Case 2 the situation is simpler. In this case one sees
immediately that $(L-(k-i0)^2)^{-1}(f\delta_T)$ is a eigenfunction
for $L$ on $D\cap\{x_2<T\}$, satisfying (1), and the proof is
complete. \vskip.1in Let $\Lambda(k)$ denote the
Dirichlet-to-Neumann operator
$$\Lambda(k)h={\partial u\over\partial x_2}\hbox{ on }x_2=T,$$
where $u$ is the outgoing solution to the boundary value problem
$Lu-k^2u=0$ in $D\cap\{x_2<T\}$, $u=h$ on $x_2=T$. Solutions to
$Lu-k^2u=0$ which vanish on $x_2=T$ and are square-integrable on
$D\cap\{0<x_1<2\pi\}\cap \{x_2<T\}$ are eigenfunctions of $L$ on
$D\cap\{x_2<T\}$ with the periodicity condition (1). When $k^2\in
S_T$, the set of eigenvalues  for $L$ on $D\cap\{x_2<T\}$, the
Dirichlet-to-Neumann map with data on $x_2=T$ is not defined.
Since the line $x_2=T$ is chosen more or less arbitrarily, for a
fixed $k$ one move $k^2$ out of $S_T$ simply by shifting $T$. The
set $S$, however, is intrinsic to the problem.

If the generalized distorted plane wave $u_+(x,k,m)$ is known for
$x_2>T$, then we know $\partial u/\partial x_2$ on $x_2=T$, and
Lemma 2 has the following corollary. \vskip.1in \noindent {\bf
Corollary 1.} The set of generalized distorted plane waves
$\{u_+(x,k,m), m\in \Bbb Z,x_2>T\}$ for fixed $k\in \Bbb R\backslash
(S_T\cup S)$, determine $\Lambda(k)$ on $x_2=T$. \vskip.1in We
want to recover $\Lambda(k)$ from the propagating modes. That
follows easily at this point. \vskip.1in \noindent{\bf Lemma 3.}
The scattering data from propagating modes,
$$\{a_m(n,k):(m+\alpha)^2<k^2 \hbox{ and }(n+\alpha)^2<k^2\}$$ given
for all $k\in \Bbb R\backslash S$, determine the distorted plane
waves in $x_2\geq T$.  \vskip.1in \noindent Proof: By (2)
$$a_m(n,k)={e^{-ix_2\sqrt{\lambda_m(k)}}\over
2\pi}\int_0^{2\pi}e^{-ix_1(m+\alpha)}v_+(x_1,T,k,n)dx_1,$$ it follows
that $a_m(n,k)$ is analytic in $k$ on the set where $v_+$ is
analytic in $k$. For fixed $m$ and $n$, $a_m(n,k)$ will be part of
the scattering data from propagating modes when $k$ is
sufficiently large. Thus for each $m$ and $n$ the scattering data
from propagating modes determine $a_m(n,k)$ on $\Bbb R\backslash
S$. Thus by (2) the propagating modes determine the generalized
distorted plane waves. \vskip.1in Combining Corollary 1 and Lemma
3, we conclude the the propagating modes determine $\Lambda(k)$
for $k\in \Bbb\backslash(S\cup S_T)$.

 \vskip.2in
\noindent \S 3. {\bf Determination of the Dirichlet-to-Neumann Map
for Wave Guides}

\vskip.1inThe arguments of the preceding section apply to the
wave guides with modifications that we give here.

Since now the unperturbed operator is $-c_0^2\Delta$, we need to
replace (3) with a representation for the outgoing fundamental solution 
for $-c_0^2\Delta$.
To obtain this representation we separate variables and use expansion 
in the
eigenfunctions (chosen to be real-valued) of the Sturm-Liouville 
problem
$$\phi_m^{\prime\prime}(x_1,k)+{k^2\over
c_0^2(x_1)}\phi_m(x_1,\mu)=\mu_m(k)\phi_m(x_1,k)$$ with
$\phi_m(0,k)=0$, $\phi_m^\prime(B,k)=0$. Using this basis and assuming 
that $k$ is chosen so that
$\mu_m(k)\neq 0, m\in {\Bbb N}$, one checks that for $f$ with bounded 
support in $x_2$
$$u(x,k)=\sum_{m=1}^\infty{1\over 2i\sqrt{\mu_m(k)}}\int_{[0,B]\times 
{\Bbb
R}}e^{i\sqrt{\mu_m(k)}|x_2-y_2|} \phi_m(x_1,k)\phi_m(y_1,k)
{f(y)\over c_0^2(y_1)} dy,\eqno{(9)}$$ is a solution to
$(L-k^2)u=f$ when $k$ is real. To see that this is the outgoing
solution we will show that $u(x,k)$ continues to a
square-integrable solution when $k$ moves into the upper half
plane. Since the boundary conditions make
$d^2/dx_1^2+k^2/c_0^2(x_1)$ self-adjoint when $k$  is real, the
functions $\phi_m(x_1,k)$ and $\mu_m(k)$ are analytic in $k$ by
Rellich's theorem. This is an elementary result here, since
$\mu_m(k)$ is a simple eigenvalue when $k$ is real. Thus for
$\epsilon>0$, if we can show that Im$\{\mu_m(k+i\epsilon)\}> 0$
when $\mu_m(k)>0$, the choice of $\sqrt{\mu_m(k+i\epsilon)}$ that
we use here (see the definitions preceding Theorem 1 in \S 1) will
make Im$\{\sqrt{\mu_m(k+i\epsilon)}>0$. However, this follows
immediately from the observation that $d\mu_m(k)/dk>0$ for $k$
real. Thus we conclude that for all $f$ for which (9) is a finite
sum, $u$ extends to a square-integrable solution to $(L-k^2)u=f$
as $k$ moves into the upper half-plane. Thus, on the complement of
the thresholds the operator $G_+(k)$,  defined by
$$G_+(k)f=
 \sum_{m=1}^\infty \phi_m(x_1,k)\int_{[0,B]\times {\Bbb
R}}{e^{i\sqrt{\mu_m(k)}|x_2-y_2|}\over 2i \sqrt{\mu_m(k)}}
\phi_m(y_1,k){f(y)\over c_0^2(y_1)}dy,$$  coincides with the limit
of $(-c_0^2\Delta -k^2I)^{-1}$ as Im$\{k\}\to 0_+$ on a dense set
of $f$. Since an easy limiting absorption argument shows that
$\lim_{\epsilon\to 0_+} (-c_0^2\Delta -k^2I)^{-1}f$  exists for
$f$ with bounded support, it follows that $G_+(k)$ is the outgoing
fundamental solution. The same construction, replacing the square
roots in (9) with their complex conjugates, leads to the incoming
fundamental solution $G_-(k)$.

As stated in \S 1, distorted plane waves for the wave guide are
obtained by solving $(L-k^2)u=0$ with the given boundary
conditions on $x_1=0$ and $x_1=B$ for $u=\Phi(x,k,m)+v_+$, where
$\Phi$ is a generalized eigenfunction for $-c_0^2(x_1)\Delta$,
i.e. $\Phi(x,m,k)=\exp(ix_2\sqrt{\mu_m(k)})\phi_m(x_1,k))$. Note
that for this to be a true distorted plane wave $\mu_m(k)$ should
be positive. However, as in \S 2 we allow \lq\lq generalized"
distorted plane waves where $\mu_m(k)<0$. As in \S 2 the
construction of outgoing distorted plane waves
$u_+(x,m,k)=\Phi(x,m,k)+v_+'(x,m,k)$ is done by limiting by
absorption. As in \S 2, $u_+$ has a representation
$u_+=\psi(x_2)\Phi + v_+$ where

$$v_+=-\lim_{\epsilon\to
0_+}(L-(k+i\epsilon)^2)^{-1}(L-k^2)(\psi(x_2)\Phi),$$ with
$L=-c^2\Delta$. Here the cutoff function $\psi\in C^\infty(\Bbb
R)$ satisfies $\psi\equiv 1$ for $|x_2|>T+1$ with support
contained in $|x_2|>T$. As before, the functions $u_+$ do not
depend on  $\psi$.  Moreover, the exceptional set S
is again the union of the thresholds and the set of $k$ for which
there is a nontrivial, square-integrable solution to $(L-k^2)u=0$
in $[0,B]\times\Bbb R$. The generalized distorted plane waves
$u_+$ have analytic continuations to Im$\{ k\}>0$, and hence as in
the case of gratings, $\{u_+(x,m,k),\ k\in I\}$, where $I$ is an
open interval in $\{k: \mu_m(k)>0\}$ determines $u_+(x,m,k)$ for
$k\in {\Bbb R}\backslash S$ (note that $\mu_m(k)>0$ for $k$
sufficiently large for each $m$). In other words the generalized
distorted plane waves are again determined by the true distorted
plane waves via analytic continuation.

The analog of (4) for wave guides is
$$\int_0^{B}\overline {f(x_1)}u_+(x_1,T,m,k)dx_1=$$
$$\int_{[0,B]\times \Bbb R}\Phi(x,m,k)
\overline{(L-k^2)\psi(L-(k-i0)^2)^{-1}(f\delta_T)}dx,$$ and this
identity shows that Corollary 1 holds for wave guides. Likewise we
have the following analog of Lemma 3. \vskip.1in \noindent {\bf
Lemma 4.} The scattering data from propagating modes,
$$\{b_m(n,k):\mu_n(k)>) \hbox{ and }\mu_m(k)>0\}$$ given
for all $k\in \Bbb R\backslash S$, determine the distorted plane
waves in $x_2\geq T$. \vskip.1in Since (2') gives,
$$b_m(n,k)=e^{-ix_2\sqrt{\mu_m(k)}}
\int_0^B\phi_m(x_1,k)v_+(x_1,T,k,n)dx_1,$$ it follows that
$a_m(n,k)$ is analytic in $k$ on the set where $u_+$ is analytic
in $k$, the proof of Lemma 3 applies here, and again conclude that
the propagating modes determine $\Lambda(k)$ for $k\in \Bbb
R\backslash(S\cup S_T)$.\vskip.2in
 \noindent \vskip.2in \noindent \S 4.
{\bf Reduction to the hyperbolic inverse problem} \vskip.2in We
will begin with the wave guide problem. Consider the hyperbolic
equation
$$v_{tt}=c^2(x)\Delta v, \ t>0 $$
in  $\{(x_1,x_2)\in [0,B]\times (-\infty,T]\}$ with zero initial
conditions, $v(x,0)=0,\ v_t(x,0)=0$, and the boundary conditions
$$v(0,x_2,t)=0,\  {\partial v\over\partial x_1}(B,x_2,t)=0, \hbox{ and 
}
v(x_1,T,t)=g(x_1,t).$$ Let $\Lambda_H$ denote the hyperbolic
Dirichlet-to-Neumann operator corresponding to this
initial-boundary value problem:
$$\Lambda_Hg={\partial v\over \partial x_2}\hbox{ for }(x_1,t)\in
[0,B]\times [0,\infty).$$ The following theorem is a particular
case of results in [B] and  [KKL] (see also [E1], [E2]).
\vskip.1in \noindent {\bf Theorem 2.} The hyperbolic
Dirichlet-to-Neumann map, $\Lambda_H$ on $x_2=T$, uniquely
determines the sound speed $c(x)$.

To deduce Theorem 1 from Theorem 2 we proceed as follows. Let
$\Lambda(k)$ be the elliptic Dirichlet-to-Neumann operator defined
previously, for $c^2(x)\Delta$, i.e. $\Lambda(k)h=\partial
u/\partial x_2$ on $x_2=T$, where $u$ is the outgoing solution to
the boundary value problem
$$c^2(x)\Delta u+k^2u=0 \hbox{ in } [0,B]\times (-\infty, T],\
u=h\hbox{ on }x_2=T$$ with the zeroDirichlet boundary condition on
$x_1=0$ and the zero Neumann condition on $x_1=B$. As we observed
earlier, $\Lambda(k)$ is analytic in $k$ off the discrete set
$S\cup S_T$. Hence, using the Fourier-Laplace transform in $t$ ,we can
recover $\Lambda_H$ from $\Lambda(k)$, given for $k_0-\epsilon
<k<k_0+\epsilon$. Since we showed in \S 3 that the propagating
modes determine $\Lambda(k)$, this completes the proof of  Theorem
1 for wave guides.

Since we have also shown for gratings  that $\Lambda(k)$ for $k\in
\Bbb R\backslash (S\cup S_T)$ is determined by scattering data
from propagating modes, the only change in the argument needed to
prove Theorem 1 for gratings is in the citation of results on the
hyperbolic Dirichlet-to-Neumann map. Here the hyperbolic
Dirichlet-to-Neumann operator, $\Lambda_H$, is defined by
$$\Lambda_Hg={\partial v\over \partial x_2}\hbox{ on } x_2=T,$$
where $v$ is the solution to
$$v_{tt}=Lv\hbox{ in } D\cap\{x_2<T\}\times \{0\leq t<\infty\}$$
satisfying the periodicity condition (1),  the initial-boundary
conditions $v(x,0)=v_t(x,0)=0$ and $v(x_1,T,t)=g$. In this setting
the uniqueness results of [B], [KKL] and [E1,2] imply that that
$\Lambda_H$ given on $x_2=T$ determine both the coefficients of
$L$, i.e. the permitivity $\epsilon(x)$, and the domain $D$. This
completes the proof of Theorem 1.

\newpage
\noindent {\bf References:}\vskip.2in

\noindent [AK] Ahluwalia, D. and J. Keller,
 Exact and asymptotic representations of the sound field in a
stratified ocean, Wave Propagation and Acoustics, Lecture Notes in
Physics, {\bf 70}(1977), 14-85.

\vskip.1in \noindent [B] Belishev, M., Boundary control in
reconstruction of manifolds and metrics (the BC method), Inverse
Problems {\bf13}(1997), R1-R45. \vskip.1in

\noindent [BDC] Bao, G., Dobson, D. and J. Cox, Mathematical
studies in rigorous grating theory, J. Optical Soc. A
{\bf12}(1995), 1029-1042. \vskip.1in

\noindent[BF] Bao, G. and A. Friedman, Inverse problems for
scattering by periodic structures, Arch. Rat. Mech. Anal. {\bf
132}(1995), 49-72. \vskip.1in

\noindent [BGWX] Buchanan, J.L., Gilbert, R.P., Mawata, C., and Y.
Xu, {\it Marine Acoustics. Direct and Inverse Problems}, SIAM,
Philadelphia, 2004. \vskip.1in

\noindent [GMX] Gilbert, R.P., Mawata, C., and Y. Xu,
Determination of a distributed inhomogeneity in a two-layered wave
guide from scattered sound, {\bf Direct and Inverse Problems of
Mathematical Physics},  R. Gilbert et al. [eds], Kluwer,
Dortrecht, 2000. \vskip.1in

\noindent[DM] Dediu, S. and J. McLaughlin, Recovering
inhomogeneities in a wave guide using eigensystem decomposition,
Inverse Problems {\bf 32}(2006), 12276-1246. \vskip.1in

\noindent [E1] Eskin, G., A new approach to the hyperbolic inverse
problems, Inverse Problems 22(2006), 815-831. \vskip.1in

\noindent [E2] Eskin, G., A new approach to the hyperbolic inverse
problems, II (Global Step), [arXiv:math AP/0701373]. \vskip.1in

 \noindent [EY] Elschner, J. and M. Yamamoto, Uniqueness
in determining polygonal, sound-hard obstacles with a single
incoming wave, Inverse Problems. 22(2006), 355-364.\vskip.1in

\noindent [H] Harari, J., Patlashenko, J. and D. Givoli,
Dirichlet-to-Neumann maps for unbounded wave guides, J. Comp.
Phys.{\bf 143}(1999), 200-223. \vskip.1in

 \noindent [HK] Hettlich, F. and A. Kirsch, Schiffer's
theorem in inverse scattering theory for periodic structures,
Inverse Problems {\bf 13}(1997), 351-361.

\vskip.1in \noindent [K] Kirsch, A., Uniqueness theorems in
inverse scattering theory for periodic structures, Inverse
Problems 10(1994), 145-152. \vskip.1in

\noindent [KKL]{} Katchalov, A., Kurylev, Y., and M. Lassas, {\it
Inverse boundary spectral problems}, Chapman \& Hall, Boca Raton,
2001. \vskip.1in

\noindent [P] R. Petit (ed.), {\it Electromagnetic Theory of
Gratings}, Topics in Current Physics {\bf 22}, Springer Verlag,
Berlin 1980.\vskip.1in

\noindent [X] Xu, Y. Inverse acoustic scattering in ocean
environments, J. Comp. Acoustics {\bf 7}(1999), 111-132.

\end